\documentclass[12pt, a4paper]{article}

\usepackage{hyperref}

\usepackage{graphicx}
\usepackage{tikz}
\usepackage{subfig}
\usepackage{pgfplots}
\usepgfplotslibrary{fillbetween}
\pgfplotsset{compat=1.18}
\usepackage{subcaption}

\usepackage[sumlimits]{amsmath}
\usepackage{amssymb,amsfonts,amsthm}
\usepackage{color,xcolor}
\usepackage{todonotes}
\usepackage{url}

\newtheorem{thm}{Theorem}
\newtheorem{cor}[thm]{Corollary}
\newtheorem{lem}[thm]{Lemma}

\newtheorem{rem}[thm]{Remark}

\def\P{\mathcal{P}}

\def\1{\mathbbm{1}}

\newcommand{\dist}[1]{\| #1\|}
\newcommand{\abs}[1]{\left|#1\right|}

\allowdisplaybreaks

\usepackage{graphicx}       

\usepackage[latin1]{inputenc}
\usepackage[T1]{fontenc}
\usepackage{microtype} 

\usepackage{cite}
\usepackage{array,colortbl}
\usepackage{amsmath,amsfonts,amssymb,bm,amsthm,mathtools,mathrsfs}
\DeclareFontFamily{U}{mathx}{\hyphenchar\font45}
\DeclareFontShape{U}{mathx}{m}{n}{
      <5> <6> <7> <8> <9> <10>
      <10.95> <12> <14.4> <17.28> <20.74> <24.88>
      mathx10
      }{}
\DeclareSymbolFont{mathx}{U}{mathx}{m}{n}
\DeclareFontSubstitution{U}{mathx}{m}{n}
\DeclareMathAccent{\widecheck}{0}{mathx}{"71}

\def\citep#1#2{\cite[{#1}]{#2}}

\usepackage{algorithm}
\usepackage{algorithmic}

\newcommand{\bbN}{{\mathbb{N}}}

\newcommand{\bbZ}{{\mathbb{Z}}}
\newcommand{\R}{{\mathbb{R}}} 
\newcommand{\Z}{{\mathbb{Z}}} 

\DeclareSymbolFont{bbold}{U}{bbold}{m}{n}
\DeclareSymbolFontAlphabet{\mathbbold}{bbold}

\newcommand{\bsx}{{\boldsymbol{x}}}
\newcommand{\bsy}{{\boldsymbol{y}}}

\newcommand{\floor}[1]{\left\lfloor #1 \right\rfloor} 

\title{Generalized Markov numbers and plane lattices}
\author{Nicolas Nagel\footnote{Corresponding author: nicolas.nagel(AT)jku.at} \ and Friedrich Pillichshammer\thanks{The authors are supported by the Austrian Science Fund (FWF) 10.55776/PAT1173425.}}
\date{}

\begin{document}

\maketitle

\begin{abstract}
We give new explicit constructions of lattice rules in the two-dimensional square and torus. It is well-known that uniformly distributed lattices (in the sense of low-discrepancy) arise from continued fractions with bounded partial denominators. We demonstrate how integers solutions of the generalized Markov equation $x^2+y^2+z^2 = kxyz-k+3$ give rise to such lattices.
\end{abstract}

\centerline{\begin{minipage}[hc]{130mm}{
}
\end{minipage}}

\section{Introduction and statement of the main result}

A \emph{Markov triple} $(x, y, z) \in \bbN^3$ is a solution in the positive integers of the \emph{Markov equation}
\begin{align}\label{eq:Markov}
    x^2+y^2+z^2 = 3xyz.
\end{align}
A \emph{Markov number} is a natural number $m \in \bbN$ which appears as a member of a Markov triple. The Markov numbers below $1000$ are $1, 2, 5, 13, 29, 34, 89, 169, 194, 233, 433, 610, 985$. The relevance of these numbers can be traced back to the works of Markov \cite{Markov}, Hurwitz \cite{Hurwitz} and Frobenius \cite{Frobenius} in their studies on Diophantine approximation and binary quadratic forms. A modern account, also describing several other appearances of the Markov numbers all throughout mathematics, is given in Aigner's book \cite{Aig13}.

In this paper, we will consider a generalization of the above equation. For any $k \in \bbN \setminus\{1,2\}$ consider the \emph{$k$-Markov equation}
\begin{align}\label{eq:k_Markov}
    x^2 + y^2 + z^2 = kxyz - k + 3.
\end{align}
We are interested in the positive integer solutions $(x, y, z) \in \bbN^3$ of this equation. 
It is straight forward to check that if $(x, y, z) \in \bbN^3$ is a solution of \eqref{eq:k_Markov}, then so is
\begin{align} \label{Vieta_involution}
    (kyz-x, y, z).
\end{align}
This is called the \emph{Vieta involution} or \emph{mutation} of $(x, y, z)$. In fact, this operation together with the symmetry of the $k$-Markov equation can be used to generate all solutions.

\begin{thm}\label{thm:k_sols}
    Let $k \in \bbN \setminus\{1,2\}$. Every positive integer solution $(x, y, z) \in \bbN^3$ of the $k$-Markov equation \eqref{eq:k_Markov} can be obtained by repeatedly applying the Vieta involution \eqref{Vieta_involution} and variable permutation to the triple $(1, 1, 1)$.
\end{thm}

The proof follows two main steps: Firstly, that all solutions $(x, y, z) \in \bbN^3$ with $y = z$ are given by $(1, 1, 1)$ and $(k-1, 1, 1)$. And secondly, that every solution $(x, y, z) \in \bbN^3$ with $x > y > z$ can be reduced (in the sense of $\max\{x, y, z\} > \max\{kyz-x, y, z\}$) by applying the Vieta involution \eqref{Vieta_involution}. We omit the details here but refer to \cite{GM22} where the details are explained for a different but similar generalization of the Markov equation \eqref{eq:Markov}.

Noting that a triple of integers is pairwise coprime if and only if the same holds for its Vieta involution \eqref{Vieta_involution}, this immediately leads to the following corollary.

\begin{cor}\label{cor:coprime}
    Let $k \in \bbN \setminus\{1,2\}$. Every positive integer solution $(x, y, z) \in \bbN^3$ of the $k$-Markov equation \eqref{eq:k_Markov} consists of pairwise coprime numbers, i.e. $\gcd(x, y) = \gcd(x, z) = \gcd(y, z) = 1$.
\end{cor}

We also note the following result on the number of solutions of \eqref{eq:k_Markov} up to a given bound.

\begin{thm}\label{thm:no_of_sols}
    Let $k \in \bbN \setminus\{1,2\}$. For $R \in \bbN$ one has
    $$
    \#\{(x, y, z) \in \{1, \dots, R\}^3: x^2+y^2+z^2 = kxyz - k + 3\} \asymp (\log R)^2.
    $$
\end{thm}

The proof of this result and further background can be found in \cite{GMR19}.

Assume now that $(N, a, b) \in \bbN^3$, $N > a \geq b$ is a solution of \eqref{eq:k_Markov}. By Corollary~\ref{cor:coprime}, there is a unique number $g \in \{1, \dots, N-1\}$ with $\gcd(g, N) = 1$ and $ag \equiv b \pmod{N}$. This number $g$ is called the \emph{characteristic value} of the triple $(N, a, b)$. Our main result is concerned with the continued fraction expansion of the quotient $g/N$.

\begin{thm}\label{thm:bound_cont_frac}
    Let $k \in \bbN\setminus\{1,2\}$. If $(N, a, b) \in \bbN^3$, $N > a \geq b$ is a solution of the $k$-Markov equation \eqref{eq:k_Markov} with characteristic value $g$, then the continued fraction expansion
    $$
    \frac gN = [0; c_1, \dots, c_n] = \frac1{\displaystyle c_1 + \frac1{\displaystyle c_2 + \frac1{\displaystyle \ddots + \frac1{\displaystyle c_n}}}}
    $$
    consists of partial denominators $c_i \in \{1,\ldots, k-1\}$ for all $i \in \{1, \dots, n\}$.
\end{thm}

The bound of $k-1$ is sharp. In fact, our methods even give more precise structural results on the partial denominators $c_i$, see Section \ref{sec:gen_Markov}.

The case of the classical Markov equation \eqref{eq:Markov} is worth considering on its own. The connection between Markov numbers and continued fraction expansions only consisting of ones and twos is well-known \cite{Aig13}, also see \cite{CS20}. For us relevant is the following, which may be compared to \cite[Theorem~C]{CS20}.

\begin{thm}\label{thm:cont_frac_3}
    Let $(N, a, b) \in \bbN^3$, $N > a \geq b$ be a solution of the Markov equation \eqref{eq:Markov} with $N \geq 5$. Let $1 \leq g < N/2$ be the unique integer with $ag \equiv \pm b \pmod{N}$. Then
    $$
    \frac gN = [0; c_1, c_2, \dots, c_n]
    $$
    with $n \in \bbN$ even, $c_i \in \{1, 2\}$, $c_1 = c_n = 2$, $c_i = c_{n+1-i}$ for all $i \in \{1, \dots, n\}$ and $c_{2i} = c_{2i+1}$ for all $i \in \{1, \dots, n/2-1\}$.
\end{thm}

Thus, the continued fraction expansion is symmetric and, except for the initial and final $2$, comes in blocks of $(1, 1)$ and $(2, 2)$. This property does not hold for the corresponding fractions $g/N$ generated from the $k$-Markov equations for $k > 3$. The proof is given in Section~\ref{sec:Markov}.

\begin{rem}\rm
To put the results in context, recall that Zaramba's conjecture from \cite[p.~76]{Zar} asks for the existence of a number $C$ with the property such that for any $N \in \mathbb{N}\setminus\{1\}$ there exists a $g \in \{1,\ldots,N-1\}$ with $\gcd(g,N)=1$ such that the continued fraction coefficients of $g/N$ are bounded by $C$. In fact, Zaremba conjectured that this holds already for $C=5$ (which would be the smallest possible choice for general $N$). So far there exist only partial answers. The conjecture is true when $N$ is a Fibonacci number (even with $C=1$). Niederreiter~\cite{Nie86} proved that Zaremba's conjecture holds true if $N$ is of the form $2^m$ or $3^m$, in which case one may even choose $C=3$, or of the form $5^m$, where $C$ can be chosen equal to 4. Niederreiter's proof provides an effective construction of $g$ with the desired properties. Bourgain and Kontorovich~\cite{boukont} and Huang~\cite{Hua15} proved that Zaremba's conjecture holds true for a set of denominators $N$ of density 1.    
\end{rem}

\paragraph{Disclosure on AI usage.} The main results (Theorems \ref{thm:bound_cont_frac} and \ref{thm:cont_frac_3}) were found and proven by the authors without AI assistance. Only the second part of the proof in Appendix~\ref{appA} was aided by the shared Austrian cloud platform Academic-AI, GPT-5.6 Sol.

\section{Applications to discrepancy theory}

Before coming to the proofs of Theorems \ref{thm:bound_cont_frac} and \ref{thm:cont_frac_3}, let us discuss some of their applications. Given a triple $(N, a, b) \in \bbN^3$, $N > a \geq b$ with $\gcd(N, a) = \gcd(N, b) = 1$, the associated two-dimensional \emph{(rank-one) lattice} is given by
\begin{align}\label{lat}
    \P_{N; a, b} \coloneqq \left\{\left(\left\{\frac{aj}N\right\}, \left\{\frac{bj}N\right\}\right): j=0, 1, \dots, N-1\right\},
\end{align}
where $\{x\} = x - \lfloor x\rfloor$ denotes the fractional part. The vector $(a,b)$ is called the \emph{generating vector} of the lattice $\P_{N; a, b}$. Different generating vectors may yield the same lattice. In fact, if $g \in \{1, \dots, N-1\}$ is such that $ag \equiv b \pmod{N}$ (and then also $\gcd(g, N) = 1$), then we may normalize $\P_{N; a, b} = \P_{N; 1, g}$, where the latter lattice has the generating vector $(1,g)$. Rank-one lattices form the basis for important classes of quasi-Monte Carlo integration rules, which are known as lattice rules. A general introduction to lattice rules can be found in the books \cite{DKP22,niesiam,SJ}.

It is desirable to construct lattices which fill out the space $[0, 1)^2$ (either interpreted as the unit square or the two-dimensional torus) as uniformly as possible. Multiple measures of uniformity have been proposed in this context, perhaps the most widely used being discrepancy. We briefly describe the main aspects.

For $\bsx = (x_1, x_2)$, $\bsy = (y_1, y_2)$ in $[0, 1)^2$ with $\bsx \leq \bsy$ componentwise, consider the box $[\bsx, \bsy) = [x_1, y_1) \times [x_2, y_2) \subseteq [0, 1)^2$. For an $N$-element point set $\P \subseteq [0, 1)^2$, its \emph{local discrepancy function} is given by
$$
\Delta_\P(\bsx, \bsy) \coloneqq \frac{\#(\P \cap [\bsx, \bsy))}N - |[\bsx, \bsy)|
$$
(where $|A|$ denotes the Lebesgue measure of a measurable $A \subseteq [0, 1)^2$) and its \emph{(extreme) discrepancy} is given by
$$
D_N(\P) \coloneqq \sup_{[\bsx, \bsy) \subseteq [0, 1)^2} |\Delta_\P(\bsx, \bsy)|.
$$
Via the \emph{Koksma--Hlawka inequality} and its generalizations, this quantity is related to the worst-case error of the quasi-Monte Carlo integration algorithm over $[0, 1)^2$, see, for example,  \cite{DKP22,DP10,LP26,kuinie,niesiam}.

It is known, due to the celebrated result of Schmidt \cite{Sch72}, that there is a universal constant $\kappa > 0$ such that
$$
D_N(\P) \geq \kappa \frac{\log N}N
$$
for all $N$-point sets $\P \subseteq [0, 1)^2$. While the points of $\P$ can be chosen in any way, the extra structure of lattices as in \eqref{lat} makes them particularly popular \cite{DKP22}. In fact, lattices attaining Schmidt's bound in an asymptotic sense can be constructed via continued fractions. Consider the lattice $\P_{N; 1, g}$ with normalized generating vector $(1,g)$ and the associated continued fraction expansion $g/N = [0; c_1, \dots, c_n]$. Setting
$$
\mathcal K\left(\frac gN\right) \coloneqq \max_{i=1, \dots, n} c_i,
$$
one obtains (see \cite[equation (5.39)]{niesiam} and the subsequent discussion)
$$
D_N(\P_{N; 1, g}) \leq \frac1N \left(\sum_{i=1}^n c_i + 1\right) = O\left(\mathcal K\left(\frac gN\right) \frac{\log N}N\right).
$$
Thus, a sequence of lattices is of optimal order for discrepancy if the partial denominators of the associated fractions $g/N$ stay bounded. From Theorem \ref{thm:bound_cont_frac} we thus obtain the optimality of lattices of the form $\P_{N; a, b}$ where $(N, a, b)$, $N > a \geq b$ is a solution of \eqref{eq:k_Markov}, see Theorem~\ref{thm:opt_disc} below.

The quantity $\mathcal K(g/N)$ may be used to bound other quantities related to the equidistribution properties of $\P_{N;1, g}$, including their Zaremba index, figure of merit or $R$-value, see \cite{DKP22,niesiam} for details.

Other notions of discrepancy, which have become more popular in recent years, include the \emph{star, extreme} and \emph{periodic $L_2$-discrepancy}. Here, the $L_\infty$-norm in the definition of $D_N(\P)$ is replaced by the $L_2$-norms
$$
L_{2, N}^{\text{star}}(\P) \coloneqq \left(\int_{\bsy \in [0, 1)^2} \Delta_\P(\boldsymbol{0}, \bsy)^2 \, \mathrm d \bsy\right)^{1/2},
$$
$$
L_{2, N}^{\text{extr}}(\P) \coloneqq \left(\iint_{\substack{\bsx, \bsy \in [0, 1)^2 \\ \bsx \leq \bsy}} \Delta_\P(\bsx, \bsy)^2 \, \mathrm d(\bsx, \bsy)\right)^{1/2}
$$
and
$$
L_{2, N}^{\text{per}}(\P) \coloneqq \left(\iint_{\bsx, \bsy \in [0, 1)^2} \Delta_\P(\bsx, \bsy)^2 \, \mathrm d(\bsx, \bsy)\right)^{1/2}
$$
respectively. For star discrepancy, one looks at the local discrepancy function of boxes $[\boldsymbol{0}, \bsy)$ anchored in $\boldsymbol{0}$. In the case of periodic $L_2$-discrepancy, one considers the local discrepancy of \emph{periodic boxes} $[\bsx, \bsy) = [x_1, y_1) \times [x_2, y_2)$ with $[x, y)$ for general $x, y \in [0, 1)$ being given possibly with wrap-around
$$
[x, y) = \begin{cases}
    [x, y) & \mbox{if } x < y, \\
    [0, y) \cup [x, 1) & \mbox{if } x \geq y.
\end{cases}
$$
As pointed out in \cite{HKP21}, periodic $L_2$-discrepancy is also equivalent to the notion of diaphony in the sense of Zinterhof~\cite{Zi}.

Similar to above, it is known that there is a universal constant $\kappa > 0$ with
$$
\min\{L_{2, N}^{\text{star}}(\P), L_{2, N}^{\text{per}}(\P)\} \geq L_{2, N}^{\text{extr}}(\P) \geq \kappa \frac{\sqrt{\log N}}N
$$
for all $N$-point sets $\P$ in $[0, 1)^2$, see \cite{HKP21, Rot54}. Furthermore, for lattices the notions of extreme and periodic $L_2$-discrepancy are equivalent and even related by the equality
$$
L_{2, N}^{\text{per}}(\P_{N; 1, g})^2 = 4 L_{2, N}^{\text{extr}}(\P_{N; 1, g})^2 + \frac{1+N^2}{18N^4},
$$
see \cite{HKP21, Nag25}. It is again known that the periodic (and thus also extreme) $L_2$-discrepancy of $\P_{N; 1, g}$ is essentially determined by the continued fraction $g/N = [0; c_1, \dots, c_n]$ via
\begin{equation}\label{perDiscCF}
L_{2, N}^{\text{per}}(\P_{N; 1, g}) \asymp \frac1N \left(\sum_{k=1}^n c_k^2\right)^{1/2}
\end{equation}
(see Appendix \ref{appA}). As above, we may use Theorem \ref{thm:bound_cont_frac} for lattices $\P_{N; a, b}$ generated from solutions of \eqref{eq:k_Markov}.

The situation for star $L_2$-discrepancy is somewhat different, where instead one has
$$
L_{2,N}^{\text{star}}(\P_{N; 1, g}) \asymp \frac1N \left(\sum_{k=1}^n c_k^2 + \left(\sum_{k=1}^n (-1)^k c_k\right)^2\right)^{1/2}
$$
see \cite{Bil13} or \cite[Proposition~8]{Bor24}. If the lattice $\P_{N; 1, g}$ is given with $N$ and $g$ as in Theorem \ref{thm:cont_frac_3}, i.e. from solutions of the classical Markov equation \eqref{eq:Markov}, we see that the alternating sum does not contribute to the asymptotic behavior of $L_{2,N}^{\text{star}}(\P_{N; 1, g})$.

\begin{figure}[t]
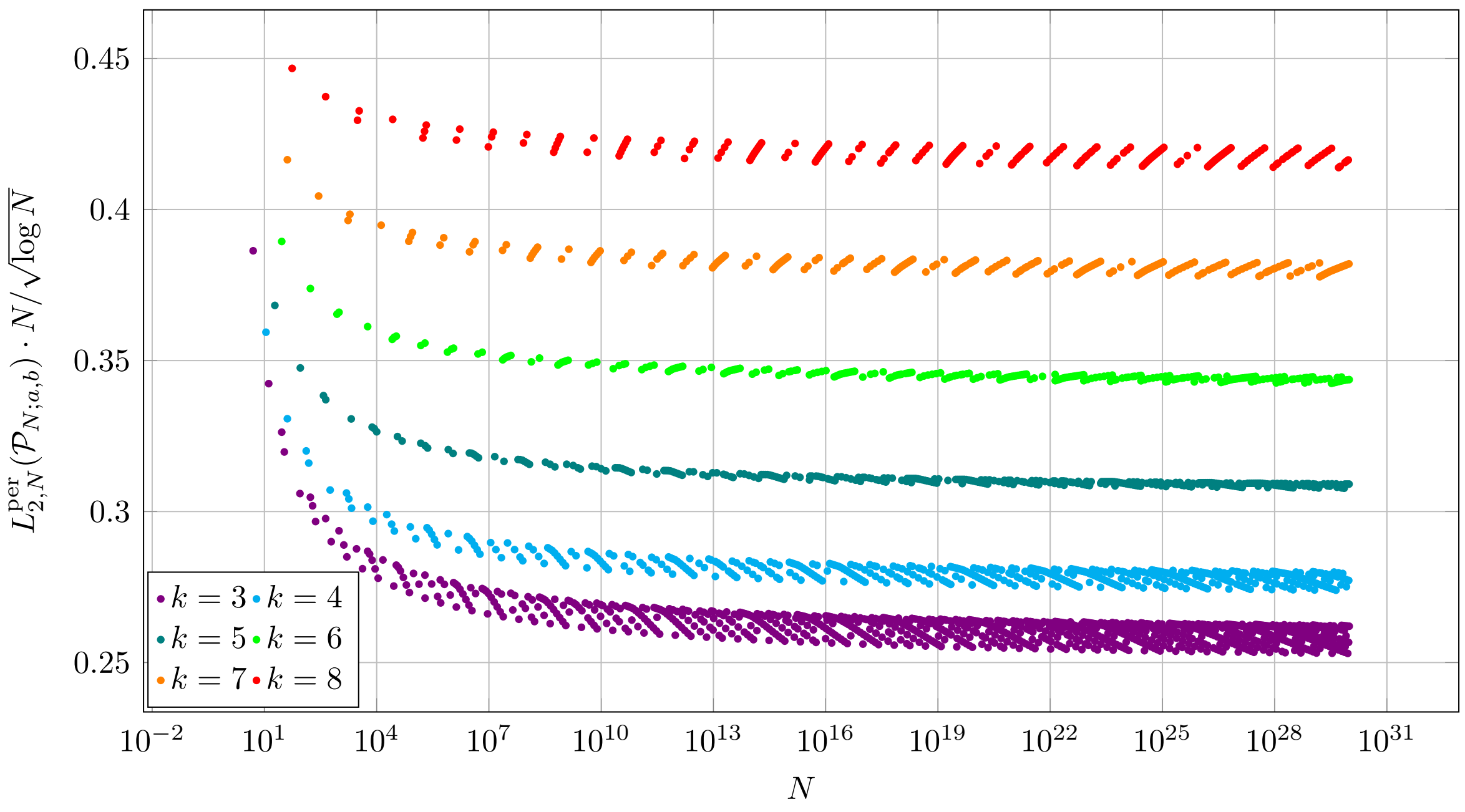


\caption{Plot of the periodic $L_2$-discrepancy normalized by a factor of $N / \sqrt{\log N}$ of the lattices $\P_{N; a, b}$ generated from the positive integer solutions $(N, a, b)$ up to $N < 10^{30}$ of the $k$-Markov equation for $k=3, \dots, 8$. The periodic $L_2$-discrepancy was calculated with an algorithm described in \cite{KN}.}
\label{fig:plots}
\end{figure}

We summarize the situation for discrepancy in the $L_\infty$- and $L_2$-sense in the following theorem. Plots of the periodic $L_2$-discrepancy for these classes of point sets are shown in Figure~\ref{fig:plots}.

\begin{thm}\label{thm:opt_disc}
    Let $k \in \bbN \setminus\{1,2\}$. If $(N, a, b) \in \bbN^3$, $N > a \geq b$ is a solution of the $k$-Markov equation \eqref{eq:k_Markov}, then
    $$
    D_N(\P_{N; a, b}) \asymp \frac{\log N}N \quad \text{and} \quad L_{2, N}^{{\rm per}}(\P_{N; a, b}) \asymp L_{2, N}^{{\rm extr}}(\P_{N; a, b}) \asymp \frac{\sqrt{\log N}}N
    $$
    with implicit constants only depending on $k$. If $k = 3$, i.e. $(N, a, b)$ is a solution of the Markov equation \eqref{eq:Markov}, one also has
    $$
    L_{2, N}^{{\rm star}}(\P_{N; a, b}) \asymp \frac{\sqrt{\log N}}N.
    $$
\end{thm}

\begin{rem}\rm
    While classical constructions of lattices via rational approximations of quadratic irrationals (including Fibonacci lattices \cite{HN}) or the method of Niederreiter \cite{Nie86} yield $O(\log R)$ many lattices of optimal order below a given bound $N \leq R$, by Theorem~\ref{thm:no_of_sols} the construction via solutions of \eqref{eq:k_Markov} yields $O((\log R)^2)$ many.
\end{rem}

\section{Tree structures}

Solutions of the Markov equation and its generalizations are most conveniently analyzed by looking at associated trees. We start with the case of the classical Markov equation \eqref{eq:Markov} (see \cite{Aig13} for further details) and then adapt this to the case of the more general $k$-Markov equation \eqref{eq:k_Markov}.

For future reference, let us state the following basic result on continued fractions.

\begin{lem}\label{lem:cont_frac_mat}
    Let $p, q \in \bbN$ with $\gcd(p, q)=1$ and let $c_1, \dots, c_n \in \bbN$ for some $n \in \bbN$. Then the continued fraction expansion $p/q = [0; c_1, \dots, c_n]$ is equivalent to the matrix factorization
    $$
\left[\begin{matrix}
    \ast & p \\
    \ast & q
\end{matrix}\right] = \left[\begin{matrix}
        0 & 1 \\
        1 & c_1
    \end{matrix}\right] \cdots \left[\begin{matrix}
        0 & 1 \\
        1 & c_n
    \end{matrix}\right].
$$
\end{lem}

\begin{proof}
Let the convergents of $[0;c_1,\ldots,c_n]$ be denoted by $p_k/q_k = [0;c_1,\ldots,c_k]$. Initializing $p_{-1}=1$, $p_0=0$, $q_{-1}=0$ and $q_0=1$, they satisfy the recurrence relations $$p_k = c_k p_{k-1}+p_{k-2} \quad \text{and} \quad q_k = c_k q_{k-1}+q_{k-2}$$ for $k \in \{1,\ldots,n\}$. In matrix notation this can be rewritten as $$ \left[\begin{matrix}
        0 & 1 \\
        1 & c_1
    \end{matrix}\right] = \left[\begin{matrix}
        p_{0} & p_1 \\
        q_{0} & q_1
    \end{matrix}\right]
$$
and then inductively
$$ \left[\begin{matrix}
        0 & 1 \\
        1 & c_1
    \end{matrix}\right] \dots \left[\begin{matrix}
        0 & 1 \\
        1 & c_{k-1}
    \end{matrix}\right] \left[\begin{matrix}
        0 & 1 \\
        1 & c_k
    \end{matrix}\right] = \left[\begin{matrix}
        p_{k-2} & p_{k-1} \\
        q_{k-2} & q_{k-1}
    \end{matrix}\right] \left[\begin{matrix}
        0 & 1 \\
        1 & c_k
    \end{matrix}\right] = \left[\begin{matrix}
        p_{k-1} & p_{k} \\
        q_{k-1} & q_{k}
    \end{matrix}\right].
$$
Thus $p/q=[0;c_1,\ldots,c_n]$ is equivalent to
$$
\left[\begin{matrix}
    \ast & p \\
    \ast & q
\end{matrix}\right] = \left[\begin{matrix}
        0 & 1 \\
        1 & c_1
    \end{matrix}\right] \cdots \left[\begin{matrix}
        0 & 1 \\
        1 & c_n
    \end{matrix}\right].
$$
where $p=p_n$, $q=q_n$ (and the starred entries are $p_{n-1}$ and $q_{n-1}$).   
\end{proof}

\subsection{Classical Markov equation $k=3$}\label{sec:Markov}

\begin{figure}[p]
\centering
\begin{tikzpicture}[scale=2.15]



\node at (0.0, 0.0) {$(1, 5, 2)$};
\node at (0.0, -0.15) {\scriptsize{$2$}};
\node at (-2.0, -0.5) {$(1, 13, 5)$};
\node at (-2.0, -0.65) {\scriptsize{$5$}};
\node at ( 2.0, -0.5) {$(5, 29, 2)$};
\node at ( 2.0, -0.65) {\scriptsize{$12$}};
\node at (-3.0, -1.0) {$(1, 34, 13)$};
\node at (-3.0, -1.15) {\scriptsize{$13$}};
\node at (-1.0, -1.0) {$(13, 194, 5)$};
\node at (-1.0, -1.15) {\scriptsize{$75$}};
\node at ( 1.0, -1.0) {$(5, 433, 29)$};
\node at ( 1.0, -1.15) {\scriptsize{$179$}};
\node at ( 3.0, -1.0) {$(29, 169, 2)$};
\node at ( 3.0, -1.15) {\scriptsize{$70$}};

\node at (-1.8, -1.5) {$(13, 7561, 194)$};
\node at (-1.8, -1.65) {\scriptsize{$2923$}};
\node at (-0.6, -1.5) {$(194, 2897, 5)$};
\node at (-0.6, -1.65) {\scriptsize{$1120$}};
\node at ( 0.6, -1.5) {$(5, 6466, 433)$};
\node at ( 0.6, -1.65) {\scriptsize{$2673$}};
\node at ( 1.8, -1.5) {$(433, 37666, 29)$};
\node at ( 1.8, -1.65) {\scriptsize{$15571$}};

\draw (-0.4, -0.1) -- (-1.6, -0.4);
\draw ( 0.4, -0.1) -- ( 1.6, -0.4);
\draw (-2.2, -0.6) -- (-2.8, -0.9);
\draw ( 2.2, -0.6) -- ( 2.8, -0.9);
\draw (-1.8, -0.6) -- (-1.2, -0.9);
\draw ( 1.8, -0.6) -- ( 1.2, -0.9);

\draw (-1.64, -1.40) -- (-1.16, -1.10);
\draw (-0.66, -1.40) -- (-0.90, -1.10);
\draw ( 1.64, -1.40) -- ( 1.16, -1.10);
\draw ( 0.66, -1.40) -- ( 0.90, -1.10);

\draw[densely dashed] (-3.1, -1.1) -- (-3.3, -1.3);
\draw[densely dashed] (-2.9, -1.1) -- (-2.7, -1.3);
\draw[densely dashed] ( 3.1, -1.1) -- ( 3.3, -1.3);
\draw[densely dashed] ( 2.9, -1.1) -- ( 2.7, -1.3);


\node at (0.0, -2.3) {\footnotesize{$\left(\left[\begin{matrix}
    0 & 1 \\
    -1 & 3
\end{matrix}\right], \left[\begin{matrix}
    2 & 5 \\
    5 & 13
\end{matrix}\right], \left[\begin{matrix}
    1 & 2 \\
    2 & 5
\end{matrix}\right]\right)$}};
\node at (-2.2, -3.0) {\footnotesize{$\left(\left[\begin{matrix}
    0 & 1 \\
    -1 & 3
\end{matrix}\right], \left[\begin{matrix}
    5 & 13 \\
    13 & 34
\end{matrix}\right], \left[\begin{matrix}
    2 & 5 \\
    5 & 13
\end{matrix}\right]\right)$}};
\node at (2.2, -3.0) {\footnotesize{$\left(\left[\begin{matrix}
    2 & 5 \\
    5 & 13
\end{matrix}\right], \left[\begin{matrix}
    12 & 29 \\
    31 & 75
\end{matrix}\right], \left[\begin{matrix}
    1 & 2 \\
    2 & 5
\end{matrix}\right]\right)$}};
\node at (-1.4, -3.7) {\footnotesize{$\left(\left[\begin{matrix}
    5 & 13 \\
    13 & 34
\end{matrix}\right], \left[\begin{matrix}
    75 & 194 \\
    196 & 507
\end{matrix}\right], \left[\begin{matrix}
    2 & 5 \\
    5 & 13
\end{matrix}\right]\right)$}};
\node at (1.4, -3.7) {\footnotesize{$\left(\left[\begin{matrix}
    2 & 5 \\
    5 & 13
\end{matrix}\right], \left[\begin{matrix}
    179 & 433 \\
    463 & 1120
\end{matrix}\right], \left[\begin{matrix}
    12 & 29 \\
    31 & 75
\end{matrix}\right]\right)$}};
\node at (-2.0, -4.4) {\footnotesize{$\left(\left[\begin{matrix}
    5 & 13 \\
    13 & 34
\end{matrix}\right], \left[\begin{matrix}
    2923 & 7561 \\
    7639 & 19760
\end{matrix}\right], \left[\begin{matrix}
    75 & 194 \\
    196 & 507
\end{matrix}\right]\right)$}};
\node at (2.0, -4.4) {\footnotesize{$\left( \left[\begin{matrix}
    179 & 433 \\
    463 & 1120
\end{matrix}\right], \left[\begin{matrix}
    15571 & 37666 \\
    40276 & 97427
\end{matrix}\right], \left[\begin{matrix}
    12 & 29 \\
    31 & 75
\end{matrix}\right]\right)$}};

\draw (-1.54,-2.79) -- (-0.66,-2.51);
\draw ( 1.54,-2.79) -- ( 0.66,-2.51);
\draw (-1.96,-3.21) -- (-1.64,-3.49);
\draw ( 1.96,-3.21) -- ( 1.64,-3.49);
\draw[densely dashed] (-2.45,-3.21) -- (-2.66,-3.42);
\draw[densely dashed] ( 2.45,-3.21) -- ( 2.66,-3.42);
\draw (1.58, -3.91) -- (1.82, -4.19);
\draw (-1.58, -3.91) -- (-1.82, -4.19);
\draw[densely dashed] ( 0.69,-3.91) -- ( 0.48,-4.12);
\draw[densely dashed] (-0.69,-3.91) -- (-0.48,-4.12);


\node at (0.0, -5.2) {$(A, AB, B)$};
\node at (-2.0, -5.7) {$(A, AAB, AB)$};
\node at ( 2.0, -5.7) {$(AB, ABB, B)$};
\node at (-3.0, -6.2) {$(A, AAAB, AAB)$};
\node at (-1.0, -6.2) {$(AAB, AABAB, AB)$};
\node at ( 3.0, -6.2) {$(ABB, ABBB, B)$};
\node at ( 1.0, -6.2) {$(AB, ABABB, ABB)$};
\node at (-1.5, -6.7) {$(AAB, AABAABAB, AABAB)$};
\node at ( 1.5, -6.7) {$(ABABB, ABABBABB, ABB)$};

\draw (-0.4, -5.3) -- (-1.6, -5.6);
\draw ( 0.4, -5.3) -- ( 1.6, -5.6);
\draw (-2.2, -5.8) -- (-2.8, -6.1);
\draw ( 2.2, -5.8) -- ( 2.8, -6.1);
\draw (-1.8, -5.8) -- (-1.2, -6.1);
\draw ( 1.8, -5.8) -- ( 1.2, -6.1);

\draw (-1.4, -6.6) -- (-1.1, -6.3);
\draw ( 1.4, -6.6) -- ( 1.1, -6.3);

\draw[densely dashed] (-3.1, -6.3) -- (-3.3, -6.5);
\draw[densely dashed] (-2.9, -6.3) -- (-2.8, -6.4);
\draw[densely dashed] ( 3.1, -6.3) -- ( 3.3, -6.5);
\draw[densely dashed] ( 2.9, -6.3) -- ( 2.8, -6.4);
\draw[densely dashed] (-0.4, -6.3) -- (-0.2, -6.5);
\draw[densely dashed] ( 0.4, -6.3) -- ( 0.2, -6.5);

\end{tikzpicture}
  \caption{The Markov, Cohn and word tree, all being connected by Theorems \ref{thm:Cohn_Markov} and \ref{thm:word_Cohn}. Below each shown Markov triple $(\ell, m, r)$ we note the canonical characteristic value $u$ as in \eqref{eq:char_val_Markov}.} \label{fig:tree3}
\end{figure}
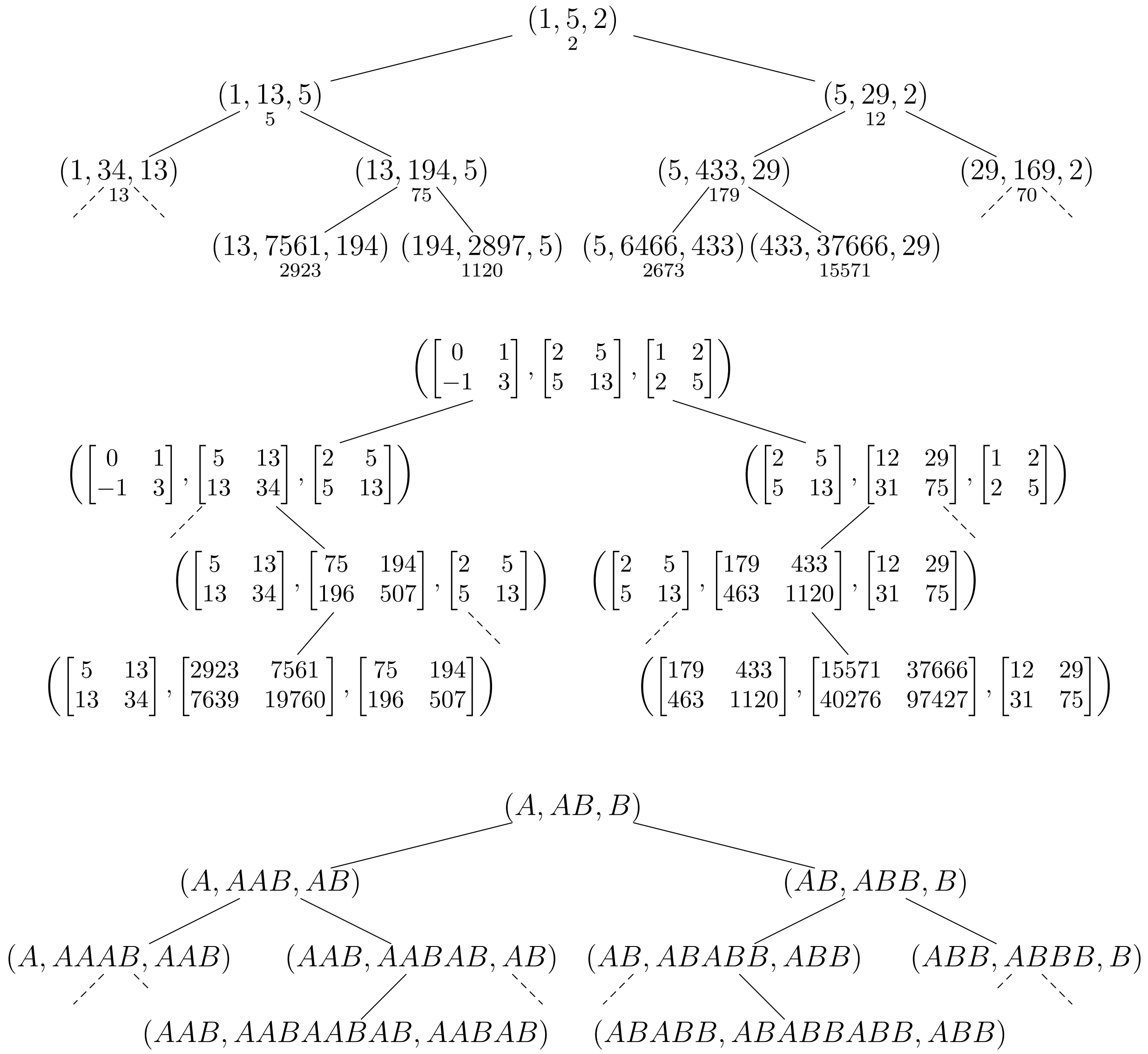

The structures that will be introduced in this section are depicted in Figure~\ref{fig:tree3}.

The \emph{Markov tree} is a rooted tree whose nodes are positive integer solution triples of the Markov equation \eqref{eq:Markov}. The root is given by $(1, 5, 2)$. For any node $(\ell, m, r)$ we define its left child to be $(\ell, 3\ell m-r, m)$ and its right child to be $(m, 3mr-\ell, r)$. By Theorem \ref{thm:k_sols}, up to permutation of the variables every solution $(x, y, z) \in \bbN^3$ of the Markov equation \eqref{eq:Markov} (except for the ``singular'' solutions $(1, 1, 1)$ and $(1, 1, 2)$) appears in this tree (and, in fact, appears exactly once in this tree, see \cite[Theorem 3.3]{Aig13}). In this way, we obtain a canonical ordering $(\ell, m, r)$ of the solutions of \eqref{eq:Markov} with the maximal element in the middle and thus also a canonical choice for the characteristic value $u \in \{1, \dots, m-1\}$ by requiring
\begin{align}\label{eq:char_val_Markov}
    \ell u \equiv r \pmod{m}.
\end{align}
In fact, one then even has $u < m/2$ (see \cite[Remark~4.14]{Aig13}), compare with the statement of Theorem~\ref{thm:cont_frac_3}.

The Markov tree is often not enough to analyze the positive integer solutions of the Markov equation \eqref{eq:Markov}. An idea going back to Cohn \cite{Coh55} is to consider an associated tree of matrices. This leads to the \emph{Cohn tree}, whose nodes consist of triples of matrices $(X, Y, Z) \in (\bbZ^{2 \times 2})^3$. The root is given by
$$
\left(\left[\begin{matrix}
    0 & 1 \\
    -1 & 3
\end{matrix}\right], \left[\begin{matrix}
    2 & 5 \\
    5 & 13
\end{matrix}\right], \left[\begin{matrix}
    1 & 2 \\
    2 & 5
\end{matrix}\right]\right)
$$
and for every node $(L, M, R)$ of this tree we define its left child to be $(L, LM, M)$ and its right child to be $(M, MR, R)$. The following result connects the Cohn tree to the Markov tree and is a direct consequence of \cite[Theorem 4.13]{Aig13}.

\begin{thm}\label{thm:Cohn_Markov}
    Let $(L, M, R)$ be a node in the Cohn tree and let $(\ell, m, r)$ be the node at the corresponding position in the Markov tree. Then $\ell, m, r$ are the top-right entries of the matrices $L, M, R$, respectively. Furthermore, the characteristic value $u$ given by \eqref{eq:char_val_Markov} is the top-left entry of $M$.
\end{thm}

Lastly, define the \emph{word tree} whose nodes consist of words over the alphabet $\{A, B\}$. The root is given by $(A, AB, B)$. If $(\lambda, \mu, \rho)$ is a node in this tree, then its left child $(\lambda, \lambda\mu, \mu)$ and right child $(\mu, \mu\rho, \rho)$ are obtained by juxtaposition. The following is part of \cite[Proposition 7.5]{Aig13}.

\begin{lem}\label{lem:palindrome}
    Let $\omega$ be an $\{A, B\}$-word appearing in some node of the word tree. Then either $\omega = A$ or $\omega = B$ or $\omega = A \sigma B$ for some (possibly empty) $\{A, B\}$-palindrome $\sigma$, i.e. an $\{A, B\}$-word which reads the same forwards as backwards.
\end{lem}

Associate any word $\omega$ appearing in the word tree to a matrix $T(\omega)$ according to
\begin{align}\label{eq:base_mat}
    T(A) \coloneqq \left[\begin{matrix}
    0 & 1 \\
    -1 & 3
\end{matrix}\right], \quad T(B) \coloneqq \left[\begin{matrix}
    1 & 2 \\
    2 & 5
\end{matrix}\right], \quad T(AB) \coloneqq \left[\begin{matrix}
    2 & 5 \\
    5 & 13
\end{matrix}\right]
\end{align}
and if $\omega = A \sigma B$ with the nonempty palindrome $\sigma = \gamma_1 \dots \gamma_\ell$, $\gamma_i \in \{A, B\}$ as in Lemma~\ref{lem:palindrome}, then
$$
T(\omega) \coloneqq \left[\begin{matrix}
        0 & 1 \\
        -1 & 3
    \end{matrix}\right] \left[\begin{matrix}
        0 & 1 \\
        1 & 2
    \end{matrix}\right] \Gamma_1 \ldots \Gamma_\ell \left[\begin{matrix}
        0 & 1 \\
        1 & 2
    \end{matrix}\right] = \left[\begin{matrix}
        1 & 2 \\
        3 & 5
    \end{matrix}\right] \Gamma_1 \ldots \Gamma_\ell \left[\begin{matrix}
        0 & 1 \\
        1 & 2
    \end{matrix}\right]
$$
where
$$
\Gamma_i = \left[\begin{matrix}
        0 & 1 \\
        1 & 1
    \end{matrix}\right]^2 = \left[\begin{matrix}
        1 & 1 \\
        1 & 2
    \end{matrix}\right] \text{ if } \gamma_i = A, \text{ and } \Gamma_i = \left[\begin{matrix}
        0 & 1 \\
        1 & 2
    \end{matrix}\right]^2 = \left[\begin{matrix}
        1 & 2 \\
        2 & 5
    \end{matrix}\right] \text{ if } \gamma_i = B.
$$
We have the following connection between the word tree and the Cohn tree.

\begin{thm}\label{thm:word_Cohn}
    Let $(\lambda, \mu, \rho)$ be a node in the word tree. Then $(T(\lambda), T(\mu), T(\rho))$ is the node at the corresponding position in the Cohn tree.
\end{thm}
\begin{proof}
    We proceed by induction on the nodes in the word tree. The relation
    $$
    (T(A), T(AB), T(B)) = \left(\left[\begin{matrix}
    0 & 1 \\
    -1 & 3
\end{matrix}\right], \left[\begin{matrix}
    2 & 5 \\
    5 & 13
\end{matrix}\right], \left[\begin{matrix}
    1 & 2 \\
    2 & 5
\end{matrix}\right]\right)
    $$
    between the roots of the trees is \eqref{eq:base_mat}. Let $(\lambda, \mu, \rho)$ be a node in the word tree and assume that $(L, M, R) = (T(\lambda), T(\mu), T(\rho))$ is the corresponding node in the Cohn tree. Let $\lambda = A \sigma_\lambda B, \mu = A \sigma_\mu B, \rho = A \sigma_\rho B$ with $\{A,B\}$-palindromes $\sigma_\lambda, \sigma_\mu, \sigma_\rho$ (the cases where $\lambda = A$ or $\rho = B$ need to be considered separately, but can be treated in an analogous manner). The left child is given by $(\lambda, \lambda\mu, \mu) = (A\sigma_\lambda B, A\sigma_\lambda BA\sigma_\mu B, A\sigma_\mu B)$ and, since
    $$
    T(\lambda\mu) = T(\lambda) \left[\begin{matrix}
        0 & 1 \\
        1 & 2
    \end{matrix}\right]^{-1} \left[\begin{matrix}
        0 & 1 \\
        1 & 2
    \end{matrix}\right]^2 \left[\begin{matrix}
        0 & 1 \\
        1 & 1
    \end{matrix}\right]^2 \left[\begin{matrix}
        1 & 2 \\
        3 & 5
    \end{matrix}\right]^{-1} T(\mu) = T(\lambda) T(\mu),
    $$
    its associated matrices $(T(\lambda), T(\lambda\mu), T(\mu))= (T(\lambda), T(\lambda) T(\mu), T(\mu))$ form the left child of $(L, M, R)$ in the Cohn tree. An analogous calculation gives the same conclusion for the right child $(\mu, \mu\rho, \rho)$, proving the claim.
\end{proof}

We are now ready to look at the continued fraction expansion of the quotients associated to a solution of the Markov equation \eqref{eq:Markov}.

\begin{proof}[Proof of Theorem \ref{thm:cont_frac_3}]
    Let $(N, a, b)$ be as in the statement, let $(\ell, m, r)$ be the corresponding node of the Markov tree with canonical characteristic value $u = g$ (see \eqref{eq:char_val_Markov}) and let $(L, M, R)$ be the corresponding node of the Cohn tree. By Theorem \ref{thm:Cohn_Markov}, Lemma \ref{lem:palindrome} and Theorem \ref{thm:word_Cohn} it holds that
    $$
    M = \left[\begin{matrix}
        u & m \\
        \ast & \ast
    \end{matrix}\right] = \left[\begin{matrix}
        0 & 1 \\
        -1 & 3
    \end{matrix}\right] \left[\begin{matrix}
        0 & 1 \\
        1 & c_1
    \end{matrix}\right] \dots \left[\begin{matrix}
        0 & 1 \\
        1 & c_n
    \end{matrix}\right],
    $$
    where $n \in \bbN$ is even and $c_i \in \{1, 2\}$ with $c_1 = c_n = 2$, $c_i = c_{n+1-i}$ for all $i \in \{1, \dots, n\}$ and $c_{2i} = c_{2i+1}$ for all $i \in \{1, \dots, n/2 - 1\}$. From this factorization and Lemma \ref{lem:cont_frac_mat} it follows that $g/N = u/m = [0; c_1, \dots, c_n]$, as desired.
\end{proof}

\subsection{Generalized Markov equation $k > 3$}\label{sec:gen_Markov}

It remains to show Theorem \ref{thm:bound_cont_frac} for $k > 3$. To do so, we extend the notions for the classical Markov equation from the previous section now to the generalized one. This needs some adjustment.

The \emph{$k$-Markov tree} is a rooted tree whose nodes consist of solution triples of the $k$-Markov equation \eqref{eq:k_Markov}. The root is given by the solution $(1, k^2-k-1, k-1)$. Any node $(\ell, m, r)$ has a left child $(\ell, k\ell m-r, m)$ and a right child $(m, kmr-\ell, r)$. Except for $(1, 1, 1)$ and $(1, k-1, 1)$, every positive integer solution of \eqref{eq:k_Markov} is represented in this tree (by Theorem \ref{thm:k_sols}). Furthermore, we obtain a natural ordering of a solution triple $(\ell, m, r)$ with its maximal element in the middle. If we use the notation $(\ell, m, r)$ we will always assume the ordering obtained from the $k$-Markov tree. The unique number $u \in \{1, \dots, m-1\}$ with
\begin{align}\label{eq:char_val}
    \ell u \equiv r \pmod{m}
\end{align}
is the canonical characteristic value of the triple $(\ell, m, r)$, compare with the paragraph before Theorem \ref{thm:bound_cont_frac}.

The \emph{$k$-matrix tree} has triples of matrices $(X, Y, Z) \in (\bbZ^{2 \times 2})^3$ as nodes. The root is given by
$$
\left(\left[\begin{matrix}
	1 & 0 \\
	0 & 1
\end{matrix}\right], \left[\begin{matrix}
1 & k-1 \\
k & k^2-k-1
\end{matrix}\right], \left[\begin{matrix}
1 & k-2 \\
1 & k-1
\end{matrix}\right]\right).
$$
If $(L, M, R)$ is a node in this tree, then its left child will be given by
$$
\left(L^\top, L^\top \left[\begin{matrix}
	0 & 1 \\
	-1 & k
\end{matrix}\right] M, \left[\begin{matrix}
	-1 & 1 \\
	0 & 1
\end{matrix}\right] M\right)
$$
and its right child by
$$
\left(M \left[\begin{matrix}
	-1 & 0 \\
	1 & 1
\end{matrix}\right], M \left[\begin{matrix}
0 & -1 \\
1 & k
\end{matrix}\right] R^\top, R^\top\right).
$$
The following lemma collects some basic properties of these triples and the proof is a simple inductive argument.

\begin{lem} \label{lem:mat_relations}
    Let
	$$
	(L, M, R) = \left(\left[\begin{matrix}
		\ell_{11} & \ell_{12} \\
		\ell_{21} & \ell_{22}
	\end{matrix}\right], \left[\begin{matrix}
		m_{11} & m_{12} \\
		m_{21} & m_{22}
	\end{matrix}\right], \left[\begin{matrix}
		r_{11} & r_{12} \\
		r_{21} & r_{22}
	\end{matrix}\right]\right)
	$$
	be a node of the $k$-matrix tree. One then has $\det L = 1$, $\det M = -1$, $\det R = 1$ and the relations
    $$
    m_{12} \ell_{22} - m_{22} \ell_{12} = r_{22} \quad \text{and} \quad m_{21} r_{22} - m_{22} r_{21} = \ell_{22}.
    $$
    Furthermore, the triple $(\ell_{22}, m_{22}, r_{22})$ is the node at the corresponding position in the $k$-Markov tree.
\end{lem}

Lastly, define the \emph{extended word tree} (not depending on the choice of $k$) whose nodes consist of triples $(\upsilon, \nu, \omega)$ of (possibly empty) words $\upsilon, \nu, \omega$ over the alphabet $\{A, B, \bar B, C\}$. The root of this tree is given by $(\varepsilon, AB, B)$, where $\varepsilon$ denotes the empty word. To define the left and right child of each node we need some auxiliary operations.
\begin{itemize}
	\item [(i)] The \emph{flip operator} $(\cdot)^\top$, which reverses the order of the letters in an $\{A, B, \bar B, C\}$-word and, in addition, turns every $B$ into $\bar B$ and vice versa. Some examples:
	$$
	(A\bar BCBBC)^\top = C\bar B\bar BCBA, \quad (\bar BAAC\bar B)^\top = BCAAB.
	$$
	
	\item [(ii)] The \emph{initial modifier} $(\cdot)^\theta$, which changes the first letter of an $\{A, B, \bar B, C\}$-word according to $A \leftrightarrow B$ and $C \leftrightarrow \bar B$. Some examples:
	$$
	(A\bar BCBBC)^\theta = B\bar BCBBC, \quad (\bar BAAC\bar B)^\theta = CAAC\bar B.
	$$
	
	\item [(iii)] The \emph{final modifier} $(\cdot)^\phi$, which changes the last letter of an $\{A, B, \bar B, C\}$-word according to $A \leftrightarrow \bar B$ and $C \leftrightarrow B$. Some examples:
	$$
	(A\bar BCBBC)^\phi = A\bar BCBBB, \quad (\bar BAAC\bar B)^\phi = \bar BAACA.
	$$
\end{itemize}
We point out that $\varepsilon = \varepsilon^\top = \varepsilon^\theta = \varepsilon^\phi$. Each of these operations is an involution. If $(\upsilon, \nu, \omega)$ is a triple in the extended word tree, then its left child is defined as $(\upsilon^\top, \upsilon^\top A \nu^\theta, \nu^\theta)$ and its right child as $(\nu^\phi, \nu^\phi A \omega^\top, \omega^\top)$. Thus, the word in the middle of such a triple is always the concatenation of the words on the left and the one on the right with an interposed $A$. The following structural result follows immediately from the recursive definition of the nodes in the tree.

\begin{lem}
    Let $(\lambda, \mu, \rho)$ be a node in the extended word tree. Then there are (possibly empty) $\{A, B, \bar B, C\}$-words $\sigma_\lambda$, $\sigma_\mu$ and $\sigma_\rho$, such that $(\lambda, \mu, \rho)$ is equal to either one of
    $$
    (A\sigma_\lambda A, A\sigma_\mu B, B\sigma_\rho B), (A\sigma_\lambda C, A\sigma_\mu \bar B, \bar B\sigma_\rho \bar B) \text{ or } (C\sigma_\lambda A, C\sigma_\mu \bar B, B\sigma_\rho \bar B).
    $$
    Here we use the convention that $\varepsilon$ is of the type $A \sigma_\lambda A$, $B$ is of the type $B\sigma_\rho B$ and $\bar B$ is of the type $\bar B\sigma_\rho \bar B$.
\end{lem}
\begin{proof}
    We proceed by induction on the nodes in the extended word tree. The root $(\varepsilon, AB, B)$ is of the form $(A\sigma_\lambda A, A\sigma_\mu B, B\sigma_\rho B)$ with the above mentioned convention. By definition, a node of the form $(A\sigma_\lambda A, A\sigma_\mu B, B\sigma_\rho B)$ has the left and right child
    $$
    (A \sigma_\lambda^\top A, A \sigma_\lambda^\top A A B \sigma_\mu B, B \sigma_\mu B) \qquad\text{ and }\qquad (A\sigma_\mu C, A \sigma_\mu C A \bar B \sigma_\rho^\top \bar B, \bar B \sigma_\rho^\top \bar B),
    $$
    respectively, both being one of the types in the statement. The corresponding compatibilities for $(A\sigma_\lambda C, A\sigma_\mu \bar B, \bar B\sigma_\rho \bar B)$ and $(C\sigma_\lambda A, C\sigma_\mu \bar B, B\sigma_\rho \bar B)$ are checked in an analogous way.
\end{proof}

Associate to every letter in $\{A, B, \bar B, C\}$ the matrices
\begin{align*}
	P_k(A) & \coloneqq \left[\begin{matrix}
		0 & 1 \\
		1 & k-1
	\end{matrix}\right], \\
	P_k(B) & \coloneqq \left[\begin{matrix}
		0 & 1 \\
		1 & 1
	\end{matrix}\right] \left[\begin{matrix}
	0 & 1 \\
	1 & k-2
	\end{matrix}\right] = \left[\begin{matrix}
	1 & k-2 \\
	1 & k-1
	\end{matrix}\right], \\
	P_k(\bar B) & \coloneqq \left[\begin{matrix}
		0 & 1 \\
		1 & k-2
	\end{matrix}\right] \left[\begin{matrix}
	0 & 1 \\
	1 & 1
	\end{matrix}\right] = \left[\begin{matrix}
	1 & 1 \\
	k-2 & k-1
	\end{matrix}\right], \\
	P_k(C) & \coloneqq \left[\begin{matrix}
		0 & 1 \\
		1 & 1
	\end{matrix}\right] \left[\begin{matrix}
	0 & 1 \\
	1 & k-3
	\end{matrix}\right] \left[\begin{matrix}
	0 & 1 \\
	1 & 1
	\end{matrix}\right] = \left[\begin{matrix}
	k-3 & k-2 \\
	k-2 & k-1
	\end{matrix}\right]
\end{align*}
and to every word $\omega = \gamma_1 \dots \gamma_\ell$, $\gamma_i \in \{A, B, \bar B, C\}$ the matrix
$$
P_k(\omega) \coloneqq P_k(\gamma_1) \dots P_k(\gamma_\ell),
$$
with $P_k(\varepsilon) = \left[\begin{matrix}
	1 & 0 \\
	0 & 1
\end{matrix}\right]$.

\begin{lem}\label{lem:Pk_ops}
	For any $\{A, B, \bar B, C\}$-word $\omega$ it holds that
    \begin{align*}
		P_k(\omega^\top) = P_k(\omega)^\top, \quad
		P_k(\omega^\theta) = \left[\begin{matrix}
			-1 & 1 \\
			0 & 1
		\end{matrix}\right] P_k(\omega), \quad
		P_k(\omega^\phi) = P_k(\omega) \left[\begin{matrix}
			-1 & 0 \\
			1 & 1
		\end{matrix}\right].
	\end{align*}
\end{lem}
\begin{proof}
    The relation for the flip operator is immediate from the definition. As for the other two relations, it suffices to note that
    $$
    P_k(A) P_k(B)^{-1} = P_k(B) P_k(A)^{-1} = P_k(C) P_k(\bar B)^{-1} = P_k(\bar B) P_k(C)^{-1} = \left[\begin{matrix}
			-1 & 1 \\
			0 & 1
		\end{matrix}\right]
    $$
    and
    $$
    P_k(A)^{-1} P_k(\bar B) = P_k(\bar B)^{-1} P_k(A) = P_k(C)^{-1} P_k(B) = P_k(B)^{-1} P_k(C) = \left[\begin{matrix}
			-1 & 0 \\
			1 & 1
		\end{matrix}\right],
    $$
    corresponding to the letter changes of the initial and final modifiers, respectively.
\end{proof}

The following result connects the extended word tree and the $k$-matrix tree.

\begin{thm} \label{thm:word_mat}
	If $(\lambda, \mu, \rho)$ is a node in the extended word tree, then $(P_k(\lambda), P_k(\mu), P_k(\rho))$ is the corresponding node in the $k$-matrix tree.
\end{thm}
\begin{proof}
    We proceed by induction on the nodes in the extended word tree. The root $(\varepsilon, AB, B)$ has the associated matrices
    \begin{align*}
        P_k(\varepsilon) & = \left[\begin{matrix}
	1 & 0 \\
	0 & 1
\end{matrix}\right], \\
P_k(AB) & = \left[\begin{matrix}
	0 & 1 \\
	1 & k-1
\end{matrix}\right]\left[\begin{matrix}
	1 & k-2 \\
	1 & k-1
\end{matrix}\right] = \left[\begin{matrix}
1 & k-1 \\
k & k^2-k-1
\end{matrix}\right], \\
P_k(B) & = \left[\begin{matrix}
1 & k-2 \\
1 & k-1
\end{matrix}\right],
    \end{align*}
    which are the matrices in the root of the $k$-matrix tree. Assume now that the statement holds for the triple $(\lambda, \mu, \rho)$ with $(L, M, R) = (P_k(\lambda), P_k(\mu), P_k(\rho))$. Using Lemma \ref{lem:Pk_ops}, for the left child in the extended word tree we calculate
    \begin{align*}
        P_k(\lambda^\top) & = L^\top, \\
        P_k(\lambda^\top A \mu^\theta) & = L^\top \left[\begin{matrix}
		0 & 1 \\
		1 & k-1
	\end{matrix}\right] \left[\begin{matrix}
		-1 & 1 \\
		0 & 1
	\end{matrix}\right] M = L^\top \left[\begin{matrix}
		0 & 1 \\
		-1 & k
	\end{matrix}\right] M, \\
    P_k(\mu^\theta) & = \left[\begin{matrix}
		-1 & 1 \\
		0 & 1
	\end{matrix}\right] M,
    \end{align*}
    giving the left child of $(L, M, R)$ in the $k$-matrix tree. In the same way we can show the compatibility with the right child, giving the claim.
\end{proof}

\begin{cor} \label{cor:mat_ineqs}
    Let $k \in \bbN \setminus \{1, 2, 3\}$ and let $(L, M, R)$ be a node in the $k$-matrix tree. If
    $$
    M = \left[\begin{matrix}
		m_{11} & m_{12} \\
		m_{21} & m_{22}
	\end{matrix}\right],
    $$
    then $1 \leq m_{11} \leq m_{12} < m_{22}$ and $1 \leq m_{11} \leq m_{21} < m_{22}$.
\end{cor}
\begin{proof}
    By Theorem \ref{thm:word_mat} there is a word $\mu = \gamma_1 \dots \gamma_{\ell}$, $\gamma_i \in \{A, B, \bar B, C\}$ with $\ell \geq 2$, such that $M = P_k(\mu)$. Then
    \begin{equation}\label{prodM}
    M = P_k(\gamma_1) \dots P_k(\gamma_{\ell}) = \left[\begin{matrix}
        0 & 1 \\
        1 & c_1
    \end{matrix}\right] \dots \left[\begin{matrix}
        0 & 1 \\
        1 & c_n
    \end{matrix}\right]
    \end{equation}
    for some $c_i \in \bbN$ (and in fact $c_i \in \{1, \dots, k-1\}$) and $n \geq \ell \ge 2$. Any matrix $M$ which is a product of at least two matrices of the type
    $$
    \left[\begin{matrix}
        0 & 1 \\
        1 & c
    \end{matrix}\right] \quad \mbox{with } c \in \bbN
    $$
    is easily seen to fulfill the inequalities in the statement.
\end{proof}

\begin{cor}\label{cor:char_val}
    Let $k \in \bbN \setminus \{1, 2, 3\}$ and let
    $$
    (L, M, R) = \left(\left[\begin{matrix}
		\ell_{11} & \ell_{12} \\
		\ell_{21} & \ell_{22}
	\end{matrix}\right], \left[\begin{matrix}
		m_{11} & m_{12} \\
		m_{21} & m_{22}
	\end{matrix}\right], \left[\begin{matrix}
		r_{11} & r_{12} \\
		r_{21} & r_{22}
	\end{matrix}\right]\right).
    $$
    Then $u = m_{12}$ is the canonical characteristic value of $(\ell_{22}, m_{22}, r_{22})$ as in \eqref{eq:char_val}.
\end{cor}
\begin{proof}
    By Lemma \ref{lem:mat_relations} we have $\ell_{22} m_{12} \equiv r_{22} \pmod{m_{22}}$ and by Corollary \ref{cor:mat_ineqs} we have $1 \leq m_{12} < m_{22}$. These are the defining conditions for the canonical characteristic value.
\end{proof}

In the same way one shows that $m_{21}$ is the unique number in $\{1, \dots, m_{22}-1\}$ with
\begin{align}\label{eq:other_char_val}
    r_{22} m_{21} \equiv \ell_{22} \pmod{m_{22}}.
\end{align}
In particular $m_{12} m_{21} \equiv 1 \pmod{m_{22}}$.

With all these results at hand, we are now in a position to analyze the continued fraction expansion of the quotient associated to a solution of the $k$-Markov equation \eqref{eq:k_Markov}.

\begin{proof}[Proof of Theorem \ref{thm:bound_cont_frac}]
    The case of $k=3$ is taken care of in Theorem \ref{thm:cont_frac_3}, hence assume from now that $k>3$. Let $(N, a, b) \in \bbN^3$, $N > a \geq b$ be a solution of the $k$-Markov equation \eqref{eq:k_Markov} with characteristic value $g$ given as before the statement of Theorem \ref{thm:bound_cont_frac}. If $a = b$, then $(k-1, 1, 1)$ (see the comment after Theorem \ref{thm:k_sols}) and $g/N = [0; k-1]$ fulfills the claim. Let now $N > a > b$. One of $(a, N, b)$ or $(b, N, a)$ is a node of the $k$-Markov tree. In the first case, $g$ is also the canonical characteristic value of $(a, N, b)$ and we thus have, by Theorem~\ref{thm:word_mat} (see also \eqref{prodM}) and Corollary~\ref{cor:char_val},
    $$
    \left[\begin{matrix}
        \ast & g \\
        \ast & N
    \end{matrix}\right] = \left[\begin{matrix}
        0 & 1 \\
        1 & c_1
    \end{matrix}\right] \dots \left[\begin{matrix}
        0 & 1 \\
        1 & c_n
    \end{matrix}\right]
    $$
    for some $c_i \in \{1, \dots, k-1\}$. Thus, using Lemma \ref{lem:cont_frac_mat}, also $g/N = [0; c_1, \dots, c_n]$ is of the desired form. If otherwise $(b, N, a)$ is a node of the $k$-Markov tree then instead
    $$
    \left[\begin{matrix}
        \ast & \ast \\
        g & N
    \end{matrix}\right] = \left[\begin{matrix}
        0 & 1 \\
        1 & c_1
    \end{matrix}\right] \dots \left[\begin{matrix}
        0 & 1 \\
        1 & c_n
    \end{matrix}\right],
    $$
    compare with \eqref{eq:other_char_val}. By taking the transpose, we get $g/N = [0; c_n, \dots, c_1]$, completing the proof.
\end{proof}

\begin{appendix}
\section{The proof of equation \eqref{perDiscCF}}\label{appA}


Although we are certain that equation \eqref{perDiscCF} is known to all experts, we have not found a proof explicitly for periodic $L_2$-discrepancy in the literature. For the sake of completeness, we provide a brief justification here. All necessary background on continued fractions can be found in \cite[Section 1.2]{Aig13}.

\begin{proof}[Proof of equation \ref{perDiscCF}]
We write $\alpha:=g/N=[0;c_1,\ldots,c_n]$. Let $p_k/q_k=[0;c_1,\ldots,c_k]$ be the $k$-th convergents of $\alpha$. By construction, $q_n=N$ and $p_n=g$. We start with the following exact formula for the periodic \(L_2\)-discrepancy for lattice point sets from \cite[Theorem~10]{HKP21}, which reads
\begin{equation*}
(N L_{2,N}^{\mathrm{per}}(\mathcal P_ {N;1,g}))^2 = \frac{1}{4q_n^2} \sum_{r=1}^{q_n-1} \frac{1}{\sin^2(\pi r/q_n) \sin^2(\pi r \alpha)} + \frac{1}{9}+\frac{1}{36q_n^2}.
\end{equation*}
For \(x\in\R\), let $\dist{x} := \min_{m\in\Z}\abs{x-m}$ denote the distance from \(x\) to the nearest integer. For \(1\leq r<q_n\), define $\delta_{q_n}(r) : = \min(r,q_n-r)$. Clearly, $\dist{r/q_n} =\delta_{q_n}(r)/q_n$. 

From the elementary inequalities $2\dist{x} \leq \abs{\sin(\pi x)} \leq \pi\dist{x}$ for $x\in\R$, it follows that $$\frac{2\delta_{q_n}(r)}{q_n} \leq \sin\left(\frac{\pi r}{q_n}\right) \leq \frac{\pi\delta_{q_n}(r)}{q_n} \quad \text{and} \quad 2\dist{r \alpha} \leq \abs{\sin(\pi r \alpha)} \leq \pi\dist{r \alpha}.$$
Consequently
$$
\frac{1}{4q_n^2} \sum_{r=1}^{q_n-1} \frac{1}{\sin^2(\pi r/q_n)\sin^2(\pi r\alpha)} \asymp \sum_{r=1}^{q_n-1} \frac{1}{\delta_{q_n}(r)^2\dist{r\alpha}^{2}},
$$
with absolute implied constants.

Since $\delta_{q_n}(q_n-r)=\delta_{q_n}(r)$ and $\dist{(q_n-r)\alpha} =  \dist{r\alpha}$, we may pair the terms \(r\) and \(q_n-r\). This way we obtain
\begin{equation*}
(N L_{2,N}^{\mathrm{per}}(\mathcal P_ {N;1,g}))^2 \asymp 
  \sum_{r=1}^{\floor{q_n/2}}
  \frac{1}{r^2\dist{r\alpha}^{2}}.
\end{equation*}
It remains to show that
\begin{equation}\label{eq:cf-sum}
\sum_{r=1}^{\floor{q_n/2}}\frac{1}{r^2\dist{r \alpha}^{2}} \asymp \sum_{k=1}^n c_k^2,
\end{equation}
where the implied constants are absolute, up to the following harmless convention for the final coefficients in the finite continued fraction expansion:
A rational number in $(0, 1)$ has two finite continued fraction expansions
$$[0;c_1,\ldots,c_n] = [0;c_1,\ldots,c_n-1,1]$$ when \(c_n\geq 2\). These two expansions change the sum of squares only by a controlled amount $c_n^2 \asymp (c_n-1)^2+1$ and thus the asymptotic estimate is independent of this convention.

The denominators of the convergents $p_k/q_k$ satisfy 
\begin{equation}\label{eq:q-recurrence}
q_k=c_kq_{k-1}+q_{k-2}
\end{equation}
(as already used in the proof of Lemma \ref{lem:cont_frac_mat}). Moreover, it is well-known that
\begin{equation}\label{eq:determinant}
p_kq_{k-1}-p_{k-1}q_k=(-1)^{k+1}.
\end{equation}
The complete-quotient identity states that $$\alpha =\frac{p_k\beta_{k+1}+p_{k-1}}{q_k\beta_{k+1}+q_{k-1}}, \quad \mbox{where } \beta_{k+1}= [c_{k+1};c_{k+2},\ldots,c_n].$$
Using \eqref{eq:determinant} we obtain
$$\| q_k\alpha\| = \frac{1}{q_k\beta_{k+1}+q_{k-1}}.$$ Since $c_{k+1}<\beta_{k+1}<c_{k+1}+1$ for $k+1 < n$, using the recurrence \eqref{eq:q-recurrence} it follows that $$\frac{1}{q_{k+1}+q_k} < \| q_k\alpha\| < \frac{1}{q_{k+1}}.$$ Consequently, 
\begin{equation}\label{asyk1}
\frac{1}{q_k^2\| q_k\alpha\|^2} \asymp \left(\frac{q_{k+1}}{q_k}\right)^2 \asymp c_{k+1}^2.
\end{equation}
This immediately shows the implied lower estimate in \eqref{eq:cf-sum}.
 
In order to prove the upper estimate, we partition the summation range of \eqref{eq:cf-sum}.
For $k\in \{1,\ldots,n\}$, define the \(k\)-th continued fraction block by
\[
  B_k:=\left\{r\in\mathbb Z:q_{k-1}\leq r<q_k\right\}.
\]
The best-approximation property of convergents states that
\begin{equation}
  \label{eq:best-approximation}
  0<d<q_k \quad \mbox{ implies } \quad
  \left\lVert d\alpha\right\rVert
  \geq
  \left\lVert q_{k-1}\alpha\right\rVert.
\end{equation}
If \(r,s\in B_k\) and \(r\neq s\), then $0<|r-s|<q_k$. It follows from \eqref{eq:best-approximation} that
$\left\lVert (r-s)\alpha\right\rVert \geq \left\lVert q_{k-1}\alpha\right\rVert$. Thus the points $\{r\alpha\}$, $r\in B_k$ are \(\left\lVert q_{k-1}\alpha\right\rVert\)-separated on the circle
\(\mathbb R/\mathbb Z\).

For \(j\in\bbN\) let now
\[
  A_{k,j}
  :=
  \left\{
    r\in B_k:
    j\left\lVert q_{k-1}\alpha\right\rVert
    \leq
    \left\lVert r\alpha\right\rVert
    <
    (j+1)\left\lVert q_{k-1}\alpha\right\rVert
  \right\}.
\]
The set
\[
  \left\{
    x\in\mathbb R/\mathbb Z:
    j\left\lVert q_{k-1}\alpha\right\rVert
    \leq
    \lVert x\rVert
    <
    (j+1)\left\lVert q_{k-1}\alpha\right\rVert
  \right\}
\]
is the union of at most two intervals, each of length
\(\left\lVert q_{k-1}\alpha\right\rVert\). Since the points $\{r\alpha\}$, $r\in B_k$
are \(\left\lVert q_{k-1}\alpha\right\rVert\)-separated, it follows that $\#A_{k,j} \leq 2$. Consequently,
\[
  \sum_{r\in B_k}
  \frac{1}{\left\lVert r\alpha\right\rVert^2}
  \leq
  \frac{2}{\left\lVert q_{k-1}\alpha\right\rVert^2}
  \sum_{j=1}^{\infty}\frac{1}{j^2}
  \lesssim
  \frac{1}{\left\lVert q_{k-1}\alpha\right\rVert^2}
\]
and since \(r\geq q_{k-1}\) for \(r\in B_k\), we thus obtain
\begin{align*}
  \sum_{q_{k-1}\leq r<q_k}
  \frac{1}{
    r^2\left\lVert r\alpha\right\rVert^2
  }
  &\leq
  \frac{1}{q_{k-1}^2}
  \sum_{r\in B_k}
  \frac{1}{\left\lVert r\alpha\right\rVert^2} \lesssim
  \frac{1}{
    q_{k-1}^2
    \left\lVert q_{k-1}\alpha\right\rVert^2
  }.
\end{align*}
Applying \eqref{asyk1}, with \(k-1\) in place of \(k\), we get
\[
  \sum_{q_{k-1}\leq r<q_k}
  \frac{1}{
    r^2\left\lVert r\alpha\right\rVert^2
  }
  \lesssim c_k^2.
\]
Summing over \(k=1,\ldots,n\) and observing that the blocks $B_k$ partition the set \(\{1,\ldots,q_n-1\}\), we obtain
\[
  \sum_{r=1}^{\lfloor q_n/2\rfloor}
  \frac{1}{
    r^2\left\lVert r\alpha\right\rVert^2
  } \leq \sum_{r=1}^{q_n-1}
  \frac{1}{
    r^2\left\lVert r\alpha\right\rVert^2
  }
  \lesssim
  \sum_{k=1}^n c_k^2,
\]
as desired.
\end{proof}

We note that number theoretic sums of the form like in \eqref{eq:cf-sum} and generalizations beyond our concrete case are also studied in \cite{beck, Bor17, HN}.

\end{appendix}

\noindent{\bf Authors' Addresses:}

\noindent Nicolas Nagel, Institut f\"{u}r Finanzmathematik und Angewandte Zahlentheorie, JKU Linz, Altenbergerstra{\ss}e 69, 4040 Linz, Austria.\\ Email: nicolas.nagel(AT)jku.at\\ ORCiD: 0009-0004-3362-3543\\

\noindent Friedrich Pillichshammer, Institut f\"{u}r Finanzmathematik und Angewandte Zahlentheorie, JKU Linz, Altenbergerstra{\ss}e 69, 4040 Linz, Austria.\\ Email: friedrich.pillichshammer(AT)jku.at\\ ORCiD: 0000-0001-6952-9218

\end{document}